\documentclass[12pt]{article}
\usepackage{amssymb,latexsym,amsmath}
\usepackage{graphicx}
\def\C{\mathbb C}
\def\R{\mathbb R}
\def\N{\mathbb N}
\def\Z{\mathbb Z}

 \newtheorem{cor}{Corollary}[section]
\newtheorem{thm}{Theorem}[section]
\newtheorem{lem}{Lemma}[section]

\newtheorem{prop}{Proposition}[section]

\date{10th April 2024}

\begin{document}
\sffamily
\title{Asymptotic values, deficiencies and Bank-Laine functions}
\author{J.K. Langley}
\maketitle

\begin{abstract}
Some results are proved concerning asymptotic and deficient values in connection with the  second order linear differential equation $y'' + Ay = 0$, in which the   coefficient $A$ is entire.
\end{abstract}

\centerline{\textit{In memory of Steve Bank, who passed away on 10th April, 1994}}

\section{Introduction}

The following notation will be used throughout. Let $A$ be an entire function,
and let $f_1, f_2$ be linearly independent solutions
of
\begin{equation}
 \label{de1}
y'' + A(z) y = 0,
\end{equation}
normalised so that their Wronskian satisfies $W(f_1, f_2) = f_1f_2' - f_1' f_2 = 1$.
Write
\begin{equation}
 \label{Udef}
 U = \frac{f_2}{f_1} , \quad
 E = f_1 f_2 = \frac{U}{U'} \, .
\end{equation}
Then $U$ is locally univalent,
and
its Schwarzian derivative
$S(U)$ satisfies
\begin{equation}
 \label{schwarzian}
 S(U) = \frac{U'''}{U'} - \frac32 \, \left( \frac{U''}{U'} \right)^2 =
 2A,
\end{equation}
while $A$ and $E$ are related by the Bank-Laine equation
\cite{BIL1}
\begin{equation}
 \label{bleq}
 4A = \left( \frac{E'}{E} \right)^2 - 2 \, \frac{E''}{E} - \frac{1}{E^2} =
 \frac{(E')^2 - 2E''E - 1}{E^2} .
\end{equation}
Because $W(f_1, f_2) = 1$, the product $E$ is a \textit{Bank-Laine function}, that is, an entire function satisfying
$E'(z) = \pm 1$ whenever $E(z) = 0$. Conversely, if $E$ is a  Bank-Laine function, or if $U$ is a locally univalent meromorphic function in $\C$,
then there exist an entire $A$ and solutions $f_1, f_2$ of (\ref{de1}) such that  (\ref{Udef}), (\ref{schwarzian})  and (\ref{bleq}) all hold \cite{BIL1,BIL3,Lai1}.

For a non-constant meromorphic function $f$ on $\C$, the terms
$$
\rho (f) = \limsup_{r \to + \infty} \frac{ \log^+ T(r, f)}{\log r} \, , \quad 
\lambda (f) = \limsup_{r \to + \infty} \frac{ \log^+ N(r, 1/f)}{\log r}  \leq \rho (f),
$$
denote its order of growth and the exponent of convergence of its zeros, while if $a \in \C$ then
$$
\delta (a, f) = \liminf_{r \to + \infty} \frac{m(r,1/(f-a))}{T(r,f)} =
1 - \limsup_{r \to + \infty} \frac{N(r,1/(f-a))}{T(r,f)}$$
is the Nevanlinna deficiency of the value $a$
\cite{Hay2}.
The following  results were proved in
\cite{BIL1,
Ros,Shen}.

\begin{thm}[\cite{BIL1,Ros,Shen}]
 \label{thmA}
With the above notation, $A$,  $E$ and $\lambda(E) =
\max\{ \lambda(f_1),\lambda(f_2) \}$
satisfy the following.\\
(i) If $A$ is a polynomial of positive degree $n $ then $\lambda(E) = \rho(E) =
(n+2)/2 $.\\
(ii) If $\lambda (E) + \rho(A) < + \infty$ then $\rho (E) < + \infty$.\\
(iii) If $\lambda(E) < \rho(A) < + \infty$
then $\rho(A) \in \N = \{ 1, 2, \ldots \}$.\\
(iv) If $A$ is transcendental and $\rho(A) \leq 1/2$ then $\lambda(E) = + \infty$,
while
\begin{equation}
 \label{roineq}
\frac1{\rho(A)} + \frac1{\lambda (E) } \leq 2 \quad \text{ if $1/2 < \rho(A) < 1$.}
\end{equation}
\end{thm}

The \textit{Bank-Laine conjecture}, to the effect
that if  $\lambda (E) + \rho(A) < + \infty$ then
$\rho (A) \in \N $, was  disproved by
 Bergweiler
and Eremenko 
\cite{bebl1,bebl2}, who constructed examples via  quasiconformal methods
for all orders $\rho(A) > 1/2$,
as well as showing
that  the inequality (\ref{roineq}) is sharp.
In a subsequent paper with Rempe \cite{BERM}, they also determined all possible orders for $A$ and $E$ when $\lambda (E) < + \infty$ and $E$ has only real zeros.

The starting point of the present paper is the following result \cite{Laozawa,Larubel}.

\begin{thm}
[\cite{Laozawa,Larubel}]
\label{thmB}
The following statements hold.\\
(I)
 If $A$ is a non-constant polynomial then there exists a path $\Gamma$ tending to infinity on which all solutions of (\ref{de1}) tend to $0$ as $z \to \infty$.\\
 (II) If $A$ is non-constant
 and $\rho(E) < + \infty$ then $0$ is an asymptotic value of
 $E$, that is, there
 exists a path
 tending to infinity on which $E(z) \to 0$ as $z \to \infty$.
\end{thm}

Statement (I) follows from
\cite[Lemma 1]{Laozawa}, but for completeness a  proof based on the Liouville transformation will be outlined in Section \ref{polynomial}.
Statement (II) was deduced from (\ref{bleq}) in \cite[Theorem 1.4]{Larubel}.

It is natural to ask whether the conclusion of (I)
holds when $A$ is transcendental. While this may seem unlikely,
it might be the case that
the conclusion of (II), which is evidently weaker than that of (I),   remains true if $\rho(E) = + \infty$.
The first four results of this paper will concern these questions and the general issue of asymptotic values of Bank-Laine functions.  The statement of these results will require the following standard terminology
\cite{Nev}.

Let $G$ be a transcendental meromorphic
function in $\C$,
and suppose that  $G(z) \to a \in \C \cup \{ \infty \}$ as $z \to \infty$ along a path
$\gamma $; then the inverse $G^{-1}$ is
said to have a transcendental singularity
over the asymptotic value $a$~\cite{BE,Nev}. If $a \in \C$ then for each $\varepsilon > 0$
there  exists a component $\Omega = \Omega( a, \varepsilon, G)$ of the
set $\{ z \in \C : |G(z) - a | < \varepsilon \}$ such that
$\gamma \setminus \Omega$ is bounded: these components are called neighbourhoods of the singularity \cite{BE}.
Two  paths $\gamma, \gamma'$ on which $G(z) \to a$ determine distinct singularities if the corresponding components
$\Omega( a, \varepsilon, G)$, $\Omega'( a, \varepsilon, G)$ are disjoint for  some $ \varepsilon > 0$.
The singularity
is called direct  \cite{BE} if $\Omega( a, \varepsilon, G)$, for some
$\varepsilon > 0$, contains finitely many zeros of $G-a$, and
indirect otherwise.

 Transcendental singularities over $\infty$ may be
classified using $1/G$ and $a=0$, and Iversen's theorem \cite{Nev} states that the inverse of every transcendental entire function has a direct singularity over $\infty$.
Finally, a direct singularity
over $\infty$ is called
logarithmic if there exists $M > 0$ such that
$w = \log G(z)$ maps
$\Omega( 0, 1/M, 1/G)$ univalently onto the half-plane ${\rm Re} \, w > \log M$.

\begin{thm}
 \label{thmEL}
 Suppose that $A$ is a transcendental entire function such that $A^{-1}$ has a logarithmic singularity
 over infinity.
 Then there exists a path $\Gamma$ tending to infinity in a neighbourhood of the singularity on which all solutions of (\ref{de1}) tend to $0$ as $z \to \infty$.
\end{thm}

In particular, Theorem \ref{thmEL} applies whenever
$A$ is a transcendental entire function in the
 Eremenko-Lyubich class $\mathcal{B}$: this   consists of those
transcendental meromorphic functions $g$ for which
the set of  all finite asymptotic
and critical
values
is bounded \cite{Ber4,EL,sixsmithEL}.
By a standard classification \cite[p.287]{Nev} of isolated singularities of the inverse
function, a transcendental singularity over $\infty$ for  
$g \in \mathcal{B}$ is automatically logarithmic.

Theorem \ref{thmEL} will be proved in Section \ref{logsing} using a version of the Liouville transformation developed in \cite{LaEL} for coefficients $A$ in the class $\mathcal{B}$.
In the absence of assumptions on the growth of $A$
and the singularities of $A^{-1}$, the following partial result will be proved in Section \ref{pfthmasym1}.

\begin{thm}
\label{thmasym1}
Let $E$ be a  Bank-Laine function associated via (\ref{bleq}) with a non-constant coefficient function $A$. Then at least one of the following holds:\\
(i)  $0$ is  a limit point of the set of asymptotic and critical
values of $E$;\\
(ii) $A$ is transcendental and $E^{-1}$ has a logarithmic singularity over $0$.
\end{thm}

 Indirect singularities of $E^{-1}$ over $0$ will be discussed in Section \ref{indirectexamples}.
Theorem  \ref{thmasym1} clearly falls far short of answering the second question raised after Theorem \ref{thmB},
namely whether every Bank-Laine function
associated with a
non-constant coefficient $A$
has $0$ as an asymptotic value.
The assumption that $A$ is non-constant is not redundant
in Theorem \ref{thmasym1}, as shown by
\begin{equation}
 \label{expBLex}
 E_1(z) = e^{az+b} + c,
 \quad a, b, c \in \C,
 \quad ac = \pm 1,
\end{equation}
for which (\ref{bleq}) gives
$A = - a^2/4$. This function has omitted values $c, \infty$, and so the  Nevanlinna deficiencies  $\delta(c, E_1)$, $ \delta(\infty, E_1)$
are each
automatically $1$
\cite{Hay2}. Furthermore, $E_1$ has no critical points, and $E_1^{-1}$ has just two singularities,   namely logarithmic singularities over $c$ and
$\infty$.
The next two theorems, to be proved in Sections \ref{pfthmunique} and \ref{pfthmsumdef2},  each
show that this Bank-Laine function is in a certain sense unique.

\begin{thm}
 \label{thmunique}
 Let $E$ be a Bank-Laine function of finite order, for which $0$ is not a limit point of the set of asymptotic and  critical
values of $E$, and assume that $E^{-1}$ has a direct transcendental singularity over some finite non-zero value $c$.
Then $A$ is constant and $E$ is given by
(\ref{expBLex}).
\end{thm}
It would be interesting to know  whether the assumptions of Theorem \ref{thmunique} can be weakened to any significant extent, but
the example $e^{e^z}$ shows that the result
is not true for infinite order.

Next, it is well known that if $A \not \equiv 0$ is a polynomial then the quotient $U$ always has finite order and
maximal sum $\sum_{c \in \C\cup \{ \infty \}}   \delta(c, U) = 2$
of Nevanlinna deficiencies \cite{Dra1,Er1, Er2}. The following result shows that
it is comparatively rare for the product $E$  to have this property.

\begin{thm}
 \label{thmsumdef2}
 Let $E$ be a Bank-Laine function of finite order
 with $\sum_{c \in \C
  }   \delta(c, E) = 1$.
 Then either $0$ is
  the only finite deficient value of $E$, or
 $A$ is constant and
 $E$ is given by (\ref{expBLex}).
\end{thm}
An example of a Bank-Laine function $E_2$ with finite order
and $\delta(0, E_2) = 1$, but with infinitely many zeros, is furnished by
\begin{equation}
 \label{E_2def}
 E_2(z) = \frac{e^{ 2 \pi i z^2  } \sin \pi z}\pi .
\end{equation}

It seems worth observing that Theorems \ref{thmasym1},
\ref{thmunique} and \ref{thmsumdef2} can all be written in terms of the quotient $U = f_2/f_1$,  and
the condition that $
A = S(U)/2$ is constant means that  $U$  is   a M\"obius
 transformation in either $z$ or  an exponential
 $e^{az}$, $a \in \C \setminus \{ 0\}$. Since $E = U/U'$,
Theorems \ref{thmasym1},
\ref{thmunique} and \ref{thmsumdef2} can be expressed as follows.

\begin{thm}
 Let $U$ be a locally univalent meromorphic function in the plane with  Schwarzian derivative $S(U)$.
 \\
 (a) If  $U'/U \in \mathcal{B}$ then either
 $S(U)$ is constant or the inverse of $U'/U$ has a logarithmic singularity over $\infty$.\\
 (b) If $U'/U$ has finite order and $U'/U \in \mathcal{B}$, and if the inverse of $U'/U$ has a direct transcendental singularity over a finite non-zero value, then $S(U)$ is constant.\\
 (c) If $U'/U$ has finite order and $\sum_{a \in \C \cup \{ \infty \}} \delta (a, U'/U) = 2$ then
 either $S(U)$ is constant or $\delta( \infty, U'/U) =
 \delta (0, U'/U) = 1$.
\end{thm}

The remaining  results of this paper involve  transcendental singularities and deficiencies for $A$, and depend on the following observation, established in \cite{Lasing2016, LaBLreal}.

\begin{thm}[\cite{Lasing2016,LaBLreal}]
 \label{thmspeiser}
 Let $A$ be a transcendental entire function of finite order, and suppose that (\ref{de1}) has linearly independent solutions $f_1, f_2$ with $\lambda(f_1 f_2) < + \infty$. Then the quotient $V$ of any two linearly independent solutions of (\ref{de1}) belongs to the Speiser class $\mathcal{S}$, that is, $V$ has finitely many asymptotic and critical values.
\end{thm}
Here it may be assumed that $W(f_1, f_2) = 1$, and \cite[Corollary 1.1]{Lasing2016} implies that $U = f_2/f_1$ has finitely many asymptotic values,
because $U'/U =
1/f_1f_2$ has finite order by Theorem \ref{thmA}.
Thus the locally univalent function
$U$ belongs to $\mathcal{S}$,
and so does  $V$, the latter being a M\"obius transformation composed with $U$.

\begin{thm}
 \label{thmAsing}
 Let $A$ be a transcendental entire function of finite order, and assume that $A^{-1}$ has a direct transcendental singularity over $0$. Let $E = f_1 f_2$, where $f_1, f_2$ are linearly independent solutions of (\ref{de1}) with $W(f_1, f_2) = 1$.
 Then
 $\lambda(E) = + \infty$.
\end{thm}
Conditions on the global asymptotic behaviour of $A$ which ensure that $\lambda(E) = + \infty$ were given in \cite[Theorem 1]{BLL1} and elsewhere, but it seems slightly surprising that the existence of one direct singularity over $0$ should  suffice in this way.
The simple example \begin{equation}
\label{standardex}
E(z) = e^z, \quad A(z) =
- (1+e^{-2z})/4,
\end{equation}
 shows that the singularity of $A^{-1}$ does need to be over $0$.
Theorem \ref{thmAsing} follows from
\cite[Corollary~1.1 and Theorem 1.6]{Lasing2016}, but also from a slightly more general result
(Proposition 9.1)
which will be proved in Section \ref{pfthmAsing}. 
The next result involves deficiencies for~$A$.

\begin{thm}
\label{thmdefA}
 Let $A$ be a transcendental entire function of finite order with $\delta(0, A) > 0$. Then
 $\lambda (f_1f_2 )
 = + \infty $ for any linearly independent solutions $f_1, f_2$ of (\ref{de1}).
\end{thm}
Chiang \cite{chiang}
proved that if $\rho(A) < + \infty$ and $\delta (0, A) = 1$ then $\lambda(f_1f_2)$ is always $+ \infty$. If $A$ has finite order and finitely many zeros then every solution $f$ of (\ref{de1}) has $\lambda(f) = + \infty$ \cite{BL1}, but
\begin{equation}
 \label{defexample}
 f(z) = e^{e^z} ,
 \quad A(z) = - e^{2z} - e^z ,
 \quad \delta (0, A) = \frac12, \quad \delta (0, f'/f) = 1,
\end{equation}
shows that
$\delta(0, A) > 0$ does not force $\lambda (f)$ to be infinite for every solution. Furthermore, (\ref{standardex})
underlines that the deficient value in Theorem \ref{thmdefA} needs to be $0$.

\begin{cor}
\label{corA}
 Let $A$ be a transcendental entire function of finite order and suppose that (\ref{de1}) has a
 non-trivial solution $f$ with $\delta(0, f'/f) > 0$. Then
 $\lambda (f_1f_2 )
 = + \infty $ for any linearly independent solutions $f_1, f_2$ of (\ref{de1}).
\end{cor}
The example (\ref{defexample}) illustrates Corollary \ref{corA}: here any solution $g$ of (\ref{de1}) for which $f, g$ are linearly independent satisfies
$$
W(f, g) = c \neq 0, \quad
\left( \frac{g}f \right)' =
\frac{c}{f^2},
\quad
g(z) = c e^{e^z}
\int^z e^{-2 e^t }
\, dt ,\quad\lambda (g) = + \infty .
$$

Bergweiler and Eremenko \cite{bebl2}
called a Bank-Laine function $E$ \textit{special} if it satisfies
$E' = 1$ whenever $E = 0$, which is equivalent to
the existence of an entire function $B$ such that
\begin{equation}
 \label{specialde}
 E' = BE + 1.
\end{equation}

\begin{cor}
\label{corB}
 Let $B$ be a
 transcendental entire function of finite order with
 $\delta(0, B) > 0$. Then every solution $E$ of
 (\ref{specialde}) has
 $\lambda(E) = + \infty$.
\end{cor}
Taking $E(z) = e^z$, $B(z) = 1 - e^{-z}$, shows that the condition
$\delta(0, B) > 0$
is not redundant in Corollary~\ref{corB}.
Theorem \ref{thmdefA} and its two corollaries will be proved in Section \ref{pfthmdefA}.


\medskip
\noindent
\textit{Acknowledgement.} The author thanks the referee for a very careful reading of the manuscript and for several extremely helpful suggestions.

\section{Preliminaries}

\begin{lem}
[\cite{EF2}, Lemma III]
\label{EFlemma}
 Let $1 < r < R < + \infty$ and let the function $g$ be meromorphic in $|z| \leq R$.
 Let $I(r)$ be a subset of $[0, 2 \pi ]$ of Lebesgue measure $\mu(r)$. Then
$$
\frac1{2 \pi} \int_{I(r)}
\log^+ |g(re^{i \theta } )| \, d \theta
\leq
\frac{11 R \mu(r)}{R - r}
\left( 1 + \log^+
\frac1{\mu(r)} \right) T(R, g).
$$
\end{lem}
\hfill$\Box$
\vspace{.1in}

The following version of Fuchs' small arcs lemma is due to Hayman
\cite[p.721]{Hay7}.

\begin{lem}
\label{lemfuchssmallarcs}
Let the function $g$ be meromorphic in $|z| \leq R$, with $g(0) = 1$. Let $\eta_1, \eta_2 $ be positive with $\eta_1 + \eta_2 < 1$. Then there exists $F_R \subseteq [0, (1-\eta_1 ) R]$, of Lebesgue measure at least $R(1- \eta_1 - \eta_2)$, such that if $r \in F_R$ and
$G_r$ is a sub-interval of $[0, 2 \pi]$ of length $m$ then
\begin{equation}
\int_{G_r} \, \left|   \frac{\partial}{\partial \theta }
\, \left( \log |g(re^{i \theta } )|
\right) \right| \, d \theta
\leq 400 \eta_1^{-2} \eta_2^{-1} T(R, g)
m \log \frac{2 \pi e}m .
\label{smallarcs}
\end{equation}
\end{lem}
\hfill$\Box$
\vspace{.1in}

\begin{lem}[\cite{LaEL}, Lemma 2.1]
 \label{lemhil2}
Suppose that $c > 0$ and $0 < \varepsilon < \pi $, and that the function $B$ is analytic, with $|1-B(z)| \leq c |z|^{-2}$, on  $\Omega' =
\{ z \in \C : 1 \leq R \leq |z| < + \infty , | \arg z | \leq \pi - \varepsilon \}$.
Then there exist $d > 0$, depending only on  $c$ and $\varepsilon$, and
linearly independent solutions $U(z), V(z)$ of
$y'' + B(z)y = 0$
satisfying
\begin{eqnarray*}
U(z) &=& e^{-iz} ( 1 + \delta_1(z) ), \quad
U'(z) = -i e^{-iz} ( 1 + \delta_2(z) ), \quad \nonumber \\
V(z) &=& e^{iz} ( 1 + \delta_3(z) ), \quad
V'(z) = i e^{iz} ( 1 + \delta_4(z) ), \nonumber\\
W(U, V) &=& 2i, \quad
|\delta_j(z)| \leq \frac{d}{|z|},
\label{H4}
\end{eqnarray*}
on
$
\Omega'' =
\{ z \in \C :  R < |z| < + \infty , | \arg z | <  \pi - \varepsilon \}
\setminus \{ z \in \C : {\rm {\rm Re } } (z) \leq 0, |{\rm Im }  (z) | \leq R \} .
$
\end{lem}
Note that a change of variables $z \to -z$ shows that Lemma \ref{lemhil2}
still holds if  $\Omega'$ and $ \Omega''$
are replaced by their reflections across the imaginary axis.
\hfill$\Box$
\vspace{.1in}

For $a \in \C$ and $r > 0$, denote by $D(a, r)$ and $S(a, r)$ the open disc and circle, respectively, of centre $a$ and radius $r$.


\section{The Liouville transformation}

\subsection{The polynomial case}
\label{polynomial}

Suppose  that $A$ is a  polynomial of degree $n > 0$ in (\ref{de1}), satisfying $A(z) = a_n z^n (1+o(1))$ as $z \to \infty$. Then
there are $n+2$  critical rays given by $\arg z = \theta^*$, where $a_n e^{i (n+2) \theta^*} $ is real and positive. Apply
the Liouville transformation
\begin{equation}
 \label{ltdef}
W(Z) = A(z)^{1/4} y(z) ,
\quad Z = \int_{z_1}^z \, A(t)^{1/2} \, dt .
\end{equation}
If $R$ and $1/\varepsilon $ are large and positive then $Z \sim 2 a_n^{1/2} z^{(n+2)/2} (n+2)^{-1} $ is univalent on the sectorial region
$$
S_1 = \{ z \in \C : \,
|z| > R, \, | \arg z - \theta^*| < 2 \pi/(n+2) -
\varepsilon \} .
$$
On the image $Z(S_1)$,
(\ref{de1}) and
(\ref{ltdef}) yield  an equation for $W$ of form
\begin{equation}
 \label{p6}
 \frac{d^2W}{dZ^2}  + (1+G(Z)) W = 0,
\end{equation}
in which
$G(Z) = O( |Z|^{-2})$.
Hille's asymptotic method
\cite{Hil1,Hil2} then delivers linearly independent
solutions $e^{\pm iZ} (1+o(1))$ on
$Z(S_1)$ and hence
 solutions given by
$
u_1(z) \sim A(z)^{-1/4}
e^{-iZ} $ and
$
u_2(z) \sim A(z)^{-1/4}
e^{iZ}
$
of (\ref{de1}) on $S_1$.
Moreover,  $Z(S_1)$ contains an unbounded real
interval $I $,
on which $|e^{\pm iZ}| = 1$,
and so the factor $A(z)^{-1/4} $ in $u_1, u_2$
ensures that every solution of (\ref{de1}) tends to $0$ as $z \to \infty$ on the pre-image $\Gamma$ of $I$ in $S_1$. This confirms assertion (I) of Theorem \ref{thmB}.
\hfill$\Box$
\vspace{.1in}


\subsection{Proof of Theorem \ref{thmEL}}\label{logsing}

The proof of Theorem \ref{thmEL} will use an analogue of the Hille-Liouville method
 developed in \cite{LaEL} for the case where the inverse of the coefficient $A$
in (\ref{de1}) has a logarithmic singularity over
$\infty$.
Suppose then that $A$ is a transcendental entire function
such that, for some $M > 0$,
a component $U_M$ of the set $\{ z \in \C : |A(z)| > M \}$ is mapped univalently by  $v = \log A(z)$  onto the half-plane $H$ given by
${\rm Re} \, v > \log M$.
It may be assumed that
$0 \not \in U_M$, and also that
$M = 1$, because
if $f$ solves (\ref{de1}) then
$g(z) = f(M^{-1/2}z)$
solves
$$
y'' + B(z) y = 0, \quad B(z) = M^{-1} A(M^{-1/2}z  ).
$$

Let $z = \phi(v)$ be the inverse mapping from
$H= \{ v \in \C : \text{Re} \, v >  0 \}$ to
$U_M = U_1
\subseteq \C \setminus \{ 0 \}$.
Bieberbach's theorem
\cite[Theorem 1.1]{Hay9}
yields, for $u \in H$,
\begin{equation}
 \label{h3}
\left| \frac{\phi''(u)}{\phi'(u)} \right| \leq \frac4{{\rm Re} \, u}.
\end{equation}


\begin{lem}
 \label{lemevphi}
 If
 $|v|$ is
 sufficiently large and $| \arg v | < \pi/3$ then
 $$
 \frac{{\rm Re} \, v}4
 \leq
 \log | e^{v/2} \phi'(v) |
 \leq
 \frac{3 \, {\rm Re} \, v}4.
 $$
\end{lem}
\textit{Proof.}
Join $1$
to $v$ by a straight line segment $M_v$,  parametrised with respect to $s = {\rm Re} \, u$. Then
the arc length estimate $ |du| \leq 3 \, ds $ holds on $M_v$,
and integration of (\ref{h3})
delivers
$$
\left| \log \frac{\phi'(v)}{\phi'(1)} \right|
\leq
\int_{M_v} \frac4{{\rm Re} \, u} \, |du|
\leq
12 \int_1^{{\rm Re}\, v}
\frac{ds}s = 12 \log ({\rm Re}\, v)
= o( {\rm Re} \, v ).
$$
\hfill$\Box$
\vspace{.1in}

Under the change of variables $v = \log A(z)$
a solution $y(z)$ of (\ref{de1}) on $U_1$ transforms to a solution $w(v) = y(z)$ on $H$ of
\begin{equation}
  \label{p2}
  w''(v) - \frac{\phi''(v)}{\phi'(v)} w'(v) + e^v \phi'(v)^2 w(v) = 0,\quad
  A(z) = e^v, \quad z = \phi(v)
  .
 \end{equation} 
For $v, v_1 \in H$,
and $z = \phi(v)$, $z_1 = \phi(v_1)$ in $U_1$,
define $Z = Z(v, v_1)$ as in
the second formula of (\ref{ltdef}) by
\begin{equation}
 \label{h1}
 Z(v, v_1) = \int_{z_1}^z A(t)^{1/2} \, dt = \int_{v_1}^v e^{u/2} \phi'(u) \, du.
\end{equation}


The following  is
Lemma 3.1 of \cite{LaEL}, proved for an arbitrary
analytic univalent $\phi: H \to
\C \setminus \{ 0 \}$.

\begin{lem}
[\cite{LaEL}]
 \label{lemfirstest}
Let $\varepsilon $  be a small positive real number. Then there exists a large positive real number $N_0$, depending  on $\varepsilon$ but not on $\phi$,
with the following property.

Let $v_0 \in H$ be such that  $S_0 = {\rm Re} \, v_0 \geq N_0$,
and define $v_1, v_2, v_3, K_2$ and $ K_3$ by 
\begin{equation}
\label{vjdef}
v_j = \frac{2^j S_0}{128} + i T_0, 
 \quad T_0 = {\rm Im} \, v_0, \quad 
K_j = \left\{ v_j + r e^{i \theta} : \, r \geq 0, \, - \frac{\pi}{2^j} \leq \theta \leq  \frac{\pi}{2^j} \right\}.
\end{equation}
Then the following three conclusions all hold:\\
(i) the function $Z = Z(v, v_1)$ in (\ref{h1}) satisfies,  for $v \in K_2$,
\begin{equation}
 \label{h2}
 Z =
Z(v, v_1 ) = \int_{v_1}^v e^{u/2} \phi'(u) \, du  = 2 e^{v/2} \phi'(v) (1 + \delta (v) ), \quad | \delta (v) | < \varepsilon ;
\end{equation}
(ii) $\psi = \psi (v, v_1) = \log Z(v, v_1) $
is univalent on a domain containing $K_3$;\\
(iii) $Z$ maps a domain $D$, with $v_0 \in D \subseteq K_3$,  univalently  onto
a region $M_3$ satisfying
\begin{equation}
 \label{Omegaimage}
Z_0 = Z(v_0, v_1) \in  M_3 = 
\{ Z \in \C : \, |Z_0|/8 < |Z| < + \infty, \, | \arg (\eta Z) |< 3 \pi /4 \} ,
\end{equation}
where $\eta = 1$ if ${\rm Re} \, Z_0 \geq 0$ and $\eta = -1$ if ${\rm Re} \, Z_0 < 0$.
\end{lem}

\hfill$\Box$
\vspace{.1in}


 Fix a small positive $\varepsilon$, assume that
 $v_0 = S_0$ is
 real, large and positive, and apply Lemma~\ref{lemfirstest}, with $T_0 = {\rm Im} \, v_0 = 0$.
Let $v_1, v_2, v_3, K_2$ and $K_3$ be as in (\ref{vjdef}),  and define $Z$ 
by (\ref{h1}). 
By Lemmas~\ref{lemevphi} and \ref{lemfirstest}, $|Z_0|$ is large in
(\ref{Omegaimage})
and
there exist $\eta \in \{ -1, 1 \}$ and
a domain $D \subseteq K_3$, both as in conclusion (iii), so that $M_3 = Z(D)$.
The change of variables 
\begin{equation}
 \label{p5}
 y(z) =
 w(v) = e^{-v/4} W(Z) = A(z)^{-1/4}
 W(Z)
\end{equation}
transforms (\ref{p2}) on $D$ to the equation
(\ref{p6}) on $M_3$. 
Here
\begin{equation}
 \label{p7}
 G(Z) = \frac1{16 e^v \phi'(v)^2} \left( 1 + 4 \, \frac{\phi''(v)}{\phi'(v)} \right),
 \quad
 |G(Z)| \leq \frac1{|Z|^2} ,
\end{equation}
on $M_3 = Z(D)$, by  (\ref{h3}), (\ref{h2}) and the fact that  $S_0 = v_0 =
{\rm Re} \, v_0$ is large. Now apply
Lemma~\ref{lemhil2}
with  
\begin{equation*}
 \label{Omeganewimage1}
 \Omega' = 
\{Z  \in \C : \, |Z_0|/4 \leq |Z| < + \infty, \, | \arg (\eta  Z) | \leq  5\pi /8 \} \subseteq M_3 , 
\end{equation*}
and let $M_4 = \Omega''$, so that  $Z_0 = Z(v_0, v_1) \in M_4 \subseteq \Omega' \subseteq M_3$, by the choice of $\eta$.
This gives
solutions $U_1(Z), U_2(Z)$ of (\ref{p6}) on
$M_4$, 
which satisfy $W(U_1, U_2) = 2i$ and 
\begin{equation}
\label{p8}
 |U_1(Z)e^{iZ} - 1| + | U_2(Z)e^{-iZ} - 1| \leq \frac{d}{|Z|} ,
\end{equation}
in which the  positive  constant $d$ is independent of  $v_0$ and $Z_0$, by (\ref{p7}).

Let $T_1$ be large and positive, and if $\eta = 1$ set $L_1 = [T_1, + \infty)$,
while if $\eta = -1 $ set
$L_1 = (-\infty, -T_1]$.
As $Z \to \infty$ on
$L_1$,
 the solutions $U_1, U_2$ of (\ref{p6})  both remain bounded, by
 (\ref{p8}), and therefore so does every solution. By
Lemmas \ref{lemevphi} and
\ref{lemfirstest},
${\rm Re} \, v$ must tend to $+ \infty$ on the pre-image $L_2$ in $H$ of $L_1$ under
$Z$.
Hence
every solution of (\ref{p2}) tends to $0$ on $L_2$, by (\ref{p5}),
and so every solution of
(\ref{de1}) tends to $0$ on
$\Gamma = \phi(L_2) \subseteq U_1$.

\hfill$\Box$
\vspace{.1in}

\subsection{An example}

The proofs in Sections \ref{polynomial} and \ref{logsing} use the same Liouville transformation $Z$ given by (\ref{ltdef}), although in the second case this is done in two stages, via
(\ref{p2}) and (\ref{p5}). In both cases the path on which all solutions of (\ref{de1}) tend to $0$ is a pre-image under $Z$
of an unbounded real interval.
For transcendental entire  $A$,
not necessarily in  $\mathcal{B}$, a local version of the Liouville transformation was developed in \cite{La5} and applied near a maximum modulus point of $A$, using estimates  from the Wiman-Valiron theory \cite{Hay5}.
But it seems worth pointing out that there exist entire  $A$ for which the method used in Sections \ref{polynomial} and \ref{logsing} cannot work, even though $Z$ is itself entire.

By a construction of Barth, Brannan and Hayman \cite{BBH}, there exists a transcendental entire function $G$
with the following property:  every unbounded connected plane set contains
a sequence $(w_n)$ tending to infinity
such that $(-1)^n
{\rm Re} \, G(w_n)
\leq |w_n|^{1/2}$.
For some $z_0 \in \C$ set
$$
A(z) = G'(z)^2, \quad
Z = H(z) = \int_{z_0}^z
A(t)^{1/2} \, dt
= G(z) - G(z_0).
$$
Suppose that
$Z$
is univalent on a domain $D$, and
maps a path $\Gamma
\subseteq D$ 
onto
a horizontal half-line $L$: it may be assumed that $\Gamma$ starts at $z_0$
and
$L = [0, + \infty)$ or
$L =
(-\infty, 0]$. Then
in either case there exist
$u_n \in \Gamma$, tending to infinity,
and
$v_n = G(u_n) $
with
$|v_n| \leq 2 |u_n|^{1/2}$ and so
$4 |u_n| \geq |v_n|^2$. Hence there must exist points
$Z_n \to \infty $  on $ L$
and $z_n \to \infty$ on $\Gamma$
such that $Z_n
= H(z_n)$ and
$$
\left( H^{-1} \right)'(Z_n) \to \infty ,
\quad
G'(z_n) \to 0,
\quad A(z_n) \to 0,
\quad A(z_n)^{-1/4}
\exp( \pm i Z_n )
\to \infty.
$$
Of course this does not mean that no path exists for this $A$ on which all solutions of (\ref{de1}) tend to $0$, only that the method of Sections \ref{polynomial} and \ref{logsing} cannot be reproduced for this case.
\hfill$\Box$
\vspace{.1in}




\section{Indirect singularities over $0$}\label{indirectexamples}

The first result of this section shows that  if $A$ is a non-constant polynomial then the transcendental singularity of $E^{-1}$
over $0$,
the existence of which is guaranteed by Theorem \ref{thmB}, is always indirect.

\begin{lem}
\label{lemindirect}
 Suppose that $E^{-1}$ has a direct transcendental singularity over $0$. Then
 $A$ is transcendental.
\end{lem}
\textit{Proof.}
Assume that $E^{-1}$ has a direct  singularity over $0$
but $A$ is a polynomial.
Then $E$ has finite order, and by standard estimates for logarithmic derivatives \cite{Gun2}
there exist $M > 0$ and a set $F_0 \subseteq [0, + \infty)$, of finite linear measure, with the property that, for $|z| \not \in F_0$,
\begin{equation}
 \label{logderest}
 \left|
 \frac{E'(z)}{E(z)} \right|
 +
\left|
 \frac{E''(z)}{E(z)} \right|
 \leq |z|^M .
 \end{equation}
Since the singularity is direct, there exist $c > 0$
and a component $V$ of
$\{ z \in \C : \, |E(z)| < c \}$
which contains no zeros of $E$, so that the formula
$$
v(z) = \log |c/E(z)|
\quad (z \in V), \quad
v(z) = 0 \quad
(z \not \in V),
$$
defines a non-constant subharmonic function on $\C$.
This $v$ is evidently bounded on a path tending to infinity, by Iversen's theorem applied to $E$, and so $v$ has lower order at least $1/2$ \cite{Hay7}.
Hence a theorem of Lewis, Rossi and Weitsman \cite{LRW},
in combination with
(\ref{bleq}) and  (\ref{logderest}),
gives a path $\gamma$ tending to infinity for  which
$$
\lim_{z \to \infty, z \in \gamma}
\frac{v(z)}{\log|z|}
= + \infty , \quad
\lim_{z \to \infty, z \in \gamma, |z| \not \in F_0}
\frac{\log |A(z)|}{\log|z|}
= + \infty
,
$$
an evident  contradiction.
\hfill$\Box$
\vspace{.1in}

\subsection{An example with $A$ transcendental}

Let $E = E_2$, where $E_2$ is given by (\ref{E_2def}).
Then $1 = \lambda (E) < \rho (E)
=2$, so
the  coefficient function $A$ is transcendental, by Theorem \ref{thmA}.
Write $z = x + iy $ with
$x, y$ real, and note first that
$e^{\pm iz} $
and $\sin \pi z$
are bounded on the strip
$-1 \leq y \leq 1$,
while
$$
{\rm Re} \, ( 2 \pi i z^2 ) = - 4 \pi xy \leq 0 \quad
\hbox{and} \quad
|  e^{2 \pi i z^2 } | \leq 1 \quad
\hbox{for $x, y \geq 0$.}
$$
In particular, $E(z) \to 0$ as $z \to \infty$ on the path $\Gamma$ given by
$z = x + i/\sqrt{x} $,
$1 \leq x < + \infty$.
Moreover, if $k \in \N$ is large
then $\sin \pi z$ is small, and so is $E(z)$,
on the arc
$z = k + iy$,
$0 \leq y \leq 1/\sqrt{k},$
joining the zero $k$ of $E$ to $\Gamma$.
Hence every neighbourhood of the singularity contains zeros of $E$ and the singularity is indirect.

 Because $E$ has finite order, the Bergweiler-Eremenko theorem \cite{BE} implies that $0$ must be a limit point of critical values of $E$, which may be verified directly as follows.
If $|z|$ is large and $E'(z) = 0$  then
 there exists $k \in \Z$ such that
$$
\sin \pi z \sim \pm \tan \pi z =
\pm \frac{i}{4z}
\sim \pm
\frac{i}{4k} , \quad z - k =
O \left( \frac1k \right)
, \quad z^2 = k^2 + O(1).
$$
This gives
$ e^{2 \pi i z^2} = O(1)$ and
$E(z) = o(1)$: thus
$0$ is the unique limit point in $\C \cup \{ \infty \}$ of the critical values of $E$.


\hfill$\Box$
\vspace{.1in}


\section{Proof of Theorem \ref{thmasym1}}
\label{pfthmasym1}

\begin{lem}
 \label{lemislands}
 Let $E$ be a Bank-Laine function and let $s > 0$ be such that $E$ has no asymptotic  or critical values $w$
 with $0 < |w| \leq 2s$.
 If $D$ is a component of
 $E^{-1}(D(0, s)) $ which contains a zero of $E$
 then $D$ is mapped univalently onto $D(0, s)$ by $E$,
 while
 $D$ has Euclidean diameter at most $24 s$ and Euclidean area at least $\pi s^2/16$.
 Furthermore, if $r$ is large then the number of such components which meet
 $D(0, r)$ is at most $64 (r/s)^2$.
\end{lem}
\textit{Proof.}
Let $D$ be such a component,
and let $D'$ be the component of $E^{-1}(D(0, 2s)) $ which contains $D$.
By \cite[p.287]{Nev},
$E$ maps $D'$ univalently
onto $D(0, 2s)$: then the inverse function $\phi $
satisfies $\phi'(0) = 1/E'(\phi(0)) = \pm 1$.
Koebe's distortion theorem
\cite[Theorem 1.3]{Hay9} gives
\begin{equation*}
 \label{koebe1}
 |\phi'(w)| \leq 12
 \quad \hbox{and} \quad
 | \phi_j(w)- \phi(0)|
\leq 12 s
\quad \hbox{for $|w| \leq s$,}
\end{equation*}
so that $D$ has diameter
at most $24 s$.
On the other hand,
Koebe's quarter theorem
\cite[Theorem~1.2]{Hay9} shows that
a disc of centre $ \phi(0)$ and radius $ s | \phi'(0)|/4
=  s/4 $ lies in
$D = \phi(D(0, s))$,
and so $D$ has Euclidean area at least $\pi s^2/16$.
Finally, if $r$ is large and $N$ of these components $D$ meet
$D(0, r)$ then they are all contained in
$D(0, 2r)$ and so
$
N  \pi s^2/16 \leq 4 \pi r^2 $.
\hfill$\Box$
\vspace{.1in}

Let $E$ and $A$ be as in the hypotheses of Theorem \ref{thmasym1}, and assume that (i) does not hold. If $0$ is an asymptotic value of $E$ then the corresponding singularity of $E^{-1}$ over $0$ is logarithmic by \cite[p.287]{Nev}, and $A$ must be transcendental by Lemma \ref{lemindirect}.

Assume henceforth that $0$ is not
an asymptotic value of $E$.
Then $E$ has infinite order, by Theorem \ref{thmB}, and there exists $s > 0$ such that $E$ has no asymptotic or
 critical values $w$ with
$|w| \leq 2s$.
Then
all components of
$E^{-1} ( D(0, 2s))$ are mapped univalently onto
$D(0, 2s)$ by $E$.

By Nevanlinna's second fundamental theorem, there exist a fixed $b \in D(0, s)$ and arbitrarily large $r$ such that $E-b$ has $N$ distinct zeros $z_1, \ldots, z_N$ in $|z| \leq r$, where
$N = n(r, b, E) \geq r^3$.
But this gives
at least $N \geq r^3$ distinct components of
$E^{-1} ( D(0, s))$, each containing a zero of $E$ and meeting $D(0, r)$, which contradicts Lemma \ref{lemislands}.

\hfill$\Box$
\vspace{.1in}


\section{Functions with maximal deficiency sum}
\label{sumdef2}

Suppose that the function
$F$ is transcendental and meromorphic in the
plane, of finite order, and that
 $\sum \delta (a, F) = 2$.
By results from
\cite{Dra1,Er1,Er2},
the order $\rho$ of $F$
satisfies
$2 \leq 2 \rho \in \N$,
and there are continuous functions
$L_1(r), L_2(r)$ such that, as $r\to\infty$,
\begin{equation}
L_1(cr) = L_1(r) ( 1 + o(1) ), \quad
L_2(cr) = L_2(r) + o(1), \quad
T(r, F) = L_1(r) r^\rho (1+o(1)) ,
\label{x4}
\end{equation}
the first two estimates holding
 uniformly for $1 \leq c \leq 2$.
Further, the results of \cite{Er1,Er2} give
asymptotic representations for $F$  as follows.

Choose a large positive constant $R_1$ and a
positive function $\eta (r)$ decreasing slowly to 0 and write, for
integer $j$ with $0 \leq j \leq 2 \rho $,
\begin{equation}
D_j = \{ z : |z| > R_1 , \quad
| \arg z - L_2(|z|) - \pi j / \rho | < \pi / 2 \rho - \eta (r) \} ,
\label{x8}
\end{equation}
so that $D_{ 2 \rho } = D_0$. Associated with each $D_j$ are
a deficient value $a_j
\in \C \cup \{ \infty \}$ of $F$ and representations for $F$
\cite{Er2}  holding outside a $C^0$ set $C_1$,
a union of discs $D(z_k, r_k)$
such that
\begin{equation}
z_k \to \infty,\quad
\sum_{|z_k| < r } r_k = o(r)
\quad \hbox{as $r \to \infty $.}
\label{x3}
\end{equation}
If  $a_j$ is finite
then, for
$z = re^{i\theta} $ in $D_j \setminus C_1$,
\begin{equation}
\log | F( r e^{ i \theta } ) - a_j | =
- \pi L_1(r) r^{ \rho }
( | \cos ( \rho ( \theta - L_2(r) )) | + o(1) ).
\label{x10}
\end{equation}
On the other hand, if $a_j = \infty$ and
$z = r e^{i\theta}$ lies in $D_j \setminus C_1$, then
(\ref{x10}) becomes
\begin{equation}
\log |F(r e^{i\theta})| = \pi L_1(r) r^{\rho}
( | \cos ( \rho ( \theta - L_2(r) ) ) | + o(1) ).
\label{x14}
\end{equation}
The following facts are  well known
\cite{Dra1,DL,Er1,Er2}
and follow from (\ref{x3}),
(\ref{x10}) and (\ref{x14}).

\begin{lem}
\label{lemsumdef2}
 The deficient values $a_j$ associated with adjacent regions $D_j$ are distinct, and if $F$ is entire then they cannot both be finite.

 Moreover,the set of $r > 1$ such that the circle $S(0, r)$ meets  $C_1$ has logarithmic density $0$, and
 if $F^{-1}$ has a direct singularity over some $b \in \C \cup\{ \infty \}$ then $b$ is one of the $a_j$.
\end{lem}

A proof of the last assertion is included here for completeness.
Suppose that $F^{-1}$ has a direct singularity over $b$, and that $b$ is not one of the $a_j$. It may be assumed that $b = \infty$, and that there exist $M > 0$ and a component $D$ of $\{ z \in \C : \, |F(z)| > M \}$ which contains no poles of $F$. This $M$ may be chosen so large that
(\ref{x10}) implies that $re^{i \theta } \not \in D$, for each $j$. A non-constant, non-negative, continuous subharmonic function $v$ may be defined by $\log |F(z)/M|$ on $D$, with $v = 0$ on the complement of $D$, and $v$ has  order at most $\rho$  since
\begin{equation}
 \label{vineq}
\max \{ v(z) : \, |z| = r \}
 \leq \frac3{2 \pi}
 \int_0^{2 \pi}
 v(2 r e^{i\theta}) \, d \theta
 \leq 3 T(2r, F) + O(1).
\end{equation}
Moreover, (\ref{x3}) and
(\ref{x10}) imply that, as $r \to + \infty$
outside a set of zero logarithmic density, the angular measure of $S(0, r) \cap D$ is $o(1)$.
A contradiction then arises from the next lemma, which is a completely standard consequence of estimates for harmonic measure from \cite[p.116]{Tsuji}.

\begin{lem}
 \label{lemtsujiest}
Let $0 < s < 1 $ and let $v$ be a continuous, non-negative, non-constant  subharmonic function on $\C$, with support $D$. Suppose that, as $r \to + \infty$
outside a set $F_0$ of zero logarithmic density,  $S(0, r) \cap D$ has angular measure at most $\pi s$. Then $v$ has order of growth at least
$1/s$.
\end{lem}
\textit{Proof.} Fix
a large $S > 0$. For $t \geq S$, let
$\theta^*(t)$ denote the angular measure of $S(0, t) \cap D$, with $\theta^*(t) = + \infty$ if
$S(0, t) \subseteq D$.
Then \cite[p.116]{Tsuji}  yields, as $r \to + \infty$,
\begin{eqnarray*}
 \log \left( \max \{ v(z) : \, |z| = 2r \} \right)
 + O(1)
 &\geq&
 \int_S^r \frac{ \pi \, dt}{t \theta^*(t)}
 \geq
 \int_{[S,r] \setminus F_0} \frac{  dt}{ts }\\
 &\geq& \frac1s \left( \log \frac{r}{S}
 -  \int_{[1, r]\cap F_0} \, \frac{dt}t \right)
\geq
\left( \frac1s - o(1)\right) \log r .
\end{eqnarray*}

\hfill$\Box$
\vspace{.1in}



\section{Proof of Theorem \ref{thmunique}}
\label{pfthmunique}

Assume that $E$ is as in the hypotheses. If $A$ is constant then $E$ clearly satisfies (\ref{expBLex}), so assume that $A$ is non-constant.
Then  $A$ is transcendental and
$E^{-1}$ has a logarithmic singularity over
$0$, by Theorem \ref{thmasym1}.
Moreover,
Lemma \ref{lemislands} yields, as $r \to + \infty$,
 $$
 n(r, 1/E) = O(r^2), \quad N(r, 1/E) = O(r^2).
 $$

\begin{lem}
 \label{lemgrowth1}
 The function $E$ satisfies
 \begin{equation}
\label{growth1}
 \liminf_{r \to + \infty}
 \frac{T(r, E)}{r^2} < + \infty.
 \end{equation}
\end{lem}
\textit{Proof.}
If (\ref{growth1})
is false then $\delta(0, E) = 1$,
which is incompatible with $E^{-1}$ having a direct singularity over $c \neq 0, \infty$, by Lemma \ref{lemsumdef2}.
\hfill$\Box$
\vspace{.1in}

Thus $E$ is an entire function of order at most $2$, with lower growth given by (\ref{growth1}),
and $E^{-1}$ has a logarithmic singularity over $0$, a direct singularity over $c \neq 0, \infty$, and at least two direct singularities over $\infty$.
This gives, for some positive $\varepsilon $ and $ M$, with  $M$ large and $\varepsilon $  so small that $E$ has no asymptotic or critical values $w$ with
$0 < |w| \leq 2 \varepsilon $:\\
(i) a component $D_0$ of the set
$\{ z \in \C : \, |E(z)| < \varepsilon \}$ mapped univalently by
$v = \log 1/E(z) $ onto
${\rm Re} \, v > \log 1/\varepsilon  $;\\
(ii) a component $D_1$ of the set
$\{ z \in \C : \, |E(z)-c| < \varepsilon \}$ containing no zeros of $E-c$;\\
(iii) two components $D_2, D_3$
of the set
$\{ z \in \C : \, |E(z)| > M \}$.
\\
Moreover, these $D_j$ are disjoint and no $D_j$
 meets any component of
$\{ z \in \C : \, |E(z)| < \varepsilon \}$ which contains a zero of $E$.

Because the singularities are direct, there exist,
for $j = 0, 1, 2, 3$, continuous,
non-constant and non-negative subharmonic functions $v_j$,  the support of $v_j$ being
$D_j$. Here
$v_0(z) = \log | \varepsilon/E(z)| $ for $z \in D_0$, with $v_0(z) = 0$ elsewhere, and
$v_1, v_2, v_3$ are defined analogously.
Moreover, (\ref{vineq}) holds for each $j$, with $v = v_j$ and $F$
replaced by $E$.

\begin{lem}

For large positive $t$ let $\theta_j(t)$ denote the angular measure of the intersection of
$D_j $ with the circle $S(0, t)$:
then
\begin{equation}
 \label{phiintest}
 \phi(t) = 2 \pi -
\sum_{j=0}^3 \theta_j(t) \geq 0,
\quad
 \int_1^{+\infty}
 \frac{\phi(t)}t \, dt < + \infty.
\end{equation}
Moreover, $E$ satisfies
\begin{equation}
 \label{lowergrowth1}
 \liminf_{r \to + \infty}
 \frac{T(r,E)}{r^2} > 0.
\end{equation}
\end{lem}
\textit{Proof.} The first assertion
of (\ref{phiintest}) holds
since the $D_j$ are disjoint. Next, a standard application of the Cauchy-Schwarz inequality gives, for large positive $t$,
say $t \geq S$,
$$
16 = \left( \sum_{j=0}^3 \theta_j(t)^{1/2}
\theta_j(t)^{-1/2} \right)^2 \leq
(2 \pi - \phi(t) )
\sum_{j=0}^3 \frac1{\theta_j(t)}
$$
and hence
$$
16 + \frac{2 \phi(t) }\pi
\leq
16 + \phi(t)
\sum_{j=0}^3 \frac1{\theta_j(t)}
\leq
2
\sum_{j=0}^3 \frac{\pi}{\theta_j(t)}.
$$
Using (\ref{growth1}),
choose a sequence
$r_n \to + \infty$ with
$T(4r_n, E) = O( r_n^2 )$.
 Since
(\ref{vineq}) holds with $v = v_j$ and $F = E$, the estimate
\cite[p.116]{Tsuji} implies that,
for large $n$,
\begin{eqnarray*}
 16 \log \frac{r_n}S +
 \frac2{\pi}
 \int_S^{r_n}
 \frac{\phi(t)}t \, dt
 &\leq&
 2 \sum_{j=0}^3
 \int_S^{r_n}
  \frac{\pi}{t \theta_j(t)} \, dt
\\
 &\leq&
 2 \sum_{j=0}^3
 \log^+ \left(
 \max \{ v_j(z) : \, |z| = 2r_n  \} \right) +O(1) \\
 &\leq&
 8 \log^+ T(4r_n, E) + O(1)
 \leq
 16 \log r_n + O(1).
\end{eqnarray*}
Since $r_n \to + \infty$ and
$\phi(t) \geq 0$,
this delivers
the second assertion of
 (\ref{phiintest}).
Furthermore, it cannot be the case that $T(4r_n, E) = o( r_n)^2$, which proves
(\ref{lowergrowth1}).
\hfill$\Box$
\vspace{.1in}

\begin{lem}
 \label{fewzeros}
 The function $E$ satisfies
 $N(r, 1/E) = o( r^2)$ as $r \to + \infty$.
\end{lem}
\textit{Proof.}
It follows from (\ref{phiintest}) that, as $r \to + \infty$,
$$
\int_1^ r t \phi(t) \, dt
\leq
\int_1^{\sqrt{r}} 2 \pi t \, dt
+ r^2 \int_{\sqrt{r}}^r
 \frac{\phi(t)}t \, dt
 \leq O(r) + o(r^2) = o(r^2).
 $$
Now suppose that $r$ is large, and that $E$ has $N$ distinct zeros $u_1, \ldots, u_N$ lying in
$|z| \leq r$. Let $C_k$ be the component of $E^{-1}(D(0, \varepsilon ))$ containing $u_k$.
Then $E$ maps $C_k$ univalently onto $D(0, \varepsilon )$,
by the choice of $\varepsilon$, and $C_j$
has Euclidean diameter at most
$24 \varepsilon$, but area at least
$\pi \varepsilon^2/16$, by Lemma~\ref{lemislands}.
Since the $C_k$ cannot meet any of the $D_j$, this implies that
$$
N \pi \varepsilon^2 /16
\leq
\int_1^ {2r} t \phi(t) \, dt + O(1) = o(r^2),
$$
which gives $n(r, 1/E) = o(r^2)$,
and $N(r, 1/E) = o(r^2)$ after integration.
\hfill$\Box$
\vspace{.1in}

But (\ref{lowergrowth1}) and Lemma \ref{fewzeros} together imply that $\delta(0, E) = 1$,
again contradicting
Lemma \ref{lemsumdef2}
and the existence
of  a direct singularity of
$E^{-1}$ over $c \neq 0, \infty$.
\hfill$\Box$
\vspace{.1in}




\section{Proof of Theorem \ref{thmsumdef2}}
\label{pfthmsumdef2}

Assume that $E$ is as in the hypotheses,
that $A$ is non-constant but that $E$ has at least one finite non-zero
deficient value
$c $. Since $E$ may be replaced by $E(cz)/c$, it may be assumed that $c=1$. Then the results of Section \ref{sumdef2} may be applied with $F = E$. Since a constant may be added to $L_2(r)$,
it may be assumed further that $a_0 = c = 1$, which forces $a_1 = \infty$,
and that the $C^0$ set $C_1$ is such that, if $z $ lies outside $C_1$ and
$|z| $ is large then, for some fixed $M > 0$,
\begin{equation}
 \label{logderest1}
 \left|
 \frac{E'(z)}{E(z)} \right|
 +
 \left|
 \frac{E'(z)}{E(z)-1} \right|
 +
\left|
 \frac{E''(z)}{E(z)} \right|
  +
\left|
 \frac{E''(z)}{E(z)-1} \right|
 +
 \left|
 \frac{A'(z)}{A(z)} \right|
  +
 \left|
 \frac{A'(z)}{A(z)+ 1/4} \right|
 \leq |z|^M .
 \end{equation}

\begin{lem}
 \label{lemrho1}
 $A$ is transcendental and
 the order $\rho$ of $E$ is $1$.
\end{lem}
\textit{Proof.}
The basic idea of the proof is as follows: in regions where $E(z)$ is very close to $1$, the coefficient $A(z)$ is close to $-1/4$, while in adjacent regions where $|E(z)|$ is large, $|A(z)|$ is not too large, in  both cases by (\ref{bleq}) and (\ref{logderest1}). The key is then to use harmonic measure to pull the first estimate for $A$
across to where $E(z)$ is large, and thereby to bound the growth of $E$.

To this end, fix a large positive $K$, and use (\ref{x3}) to determine the following:\\
(i) an unbounded set $F_1 \subseteq [1, + \infty)$
such that
for $r \in F_1$ the circles
$S(0,  K^{\pm 2 } r)$ do not meet $C_1$;\\
(ii)
for each $r \in F_1$,
a small positive $\sigma = \sigma(r)$ such that
\begin{equation}
 \label{sigmaest1}
\lim_{r \to + \infty, r \in F_1} \sigma (r) = 0,
\end{equation}
for which
the radial segments
\begin{eqnarray*}
\Omega_0 &=&
\{ z = t e^{i (L_2(r) + \sigma (r))} :
\, K^{-2} r \leq t \leq K^2 r \}, \\
\Omega_1 &=&
\{ z = t e^{i (L_2(r) + \pi /\rho + \sigma (r)) } :
\, K^{-2} r \leq t \leq K^2 r \},\\
\Omega_2 &=&
\{ z = t e^{i (L_2(r) + 3\pi /4\rho + \sigma(r)) } :
\, K^{-1} r \leq t \leq K r \},
\end{eqnarray*}
do not meet $C_1$.
Here (\ref{x4}) and (\ref{x8})  imply that  $\Omega_0$ lies within $D_0$, while
 $\Omega_1, \Omega_2$ lie in $D_1$.

 Let
$\Omega $ be the domain
whose boundary consists of
$\Omega_0, \Omega_1$ and the two shorter arcs $\Omega_3, \Omega_4$ of the circles $S(0, K^2 r),
S(0, K^{-2} r)$, respectively, joining
$\Omega_0$ to $ \Omega_1$.
In order  to estimate $A$ on
$\partial \Omega$,
consider first $z \in \Omega_0$, for large $r \in F_1$.
Then
(\ref{x4}),
(\ref{x10}) and
(\ref{sigmaest1}) deliver
$$
 \log | E( z ) - 1 | \leq
- \frac{15\pi}{16}  L_1(|z|) |z|^{ \rho }
\leq
- \frac{7\pi}{8}  L_1(r) |z|^{ \rho } ,
$$
which on combination with
(\ref{logderest1})
yields, for $j = 1,2$,
$$
\log \left| \frac{E^{(j)}( z )}
{E( z )}
\right| =
\log \left| \frac{E^{(j)}( z )}
{E( z )-1} \cdot
\frac{E(z)-1}{E(z)}
\right|
\leq
- \frac{3\pi}4  L_1(r)
|z|^{ \rho }  .
$$
In view of (\ref{bleq})
and (\ref{logderest1}), these estimates lead to
\begin{equation}
 \label{est1}
\log | A( z ) + 1/4 | \leq
- \frac{\pi}2  L_1(r) |z|^{ \rho } \quad \hbox{and}
\quad
\log |A'(z)| \leq
 - \frac{\pi}4  L_1(r) |z|^{ \rho }
\quad
\hbox{ for $z \in \Omega_0$.}
\end{equation}
Note that (\ref{est1})
shows that $A$, which is assumed to be non-constant, cannot be a polynomial.


Next, $|E(z)|$ is large on $\Omega_1$, by (\ref{x14}), ($\ref{sigmaest1}$) and the fact that $a_1=\infty$. Thus
(\ref{bleq}) and (\ref{logderest1}) yield
\begin{equation}
 \label{est2}
\log | A( z  ) | \leq
3M \log |z| \leq
3M \log K^2 r \quad
\hbox{ for $z \in \Omega_1$.}
\end{equation}

It remains to estimate $A$ on $\Omega_3, \Omega_4$. Now, (\ref{bleq}) and (\ref{x4}) give a
routine estimate, for $z \in \partial \Omega$,
\begin{eqnarray}
\log |A(z)| &\leq& \log M(K^2 r , A) \leq
\left( \frac{K+1}{K-1} \right)  T(K^3 r, A) \nonumber \\
&\leq& 2
\left( \frac{K+1}{K-1} \right)  T(K^3 r, E) + O( \log K^3 r) 
\leq 4 \left( \frac{K+1}{K-1} \right)
L_1(r) (K^3r)^\rho.
\label{est3}
\end{eqnarray}
However, on most of these arcs $\Omega_3$ and  $\Omega_4$, either $E(z)$ is close to $1$ or $|E(z)|$ is large, by (\ref{x10}) and (\ref{x14}),
and so an upper bound for $|A(z)|$ follows from (\ref{bleq}).
Thus (\ref{logderest1}) delivers
\begin{equation}
 \log |A(z)| \leq
 3M \log K^2r
 \label{est4}
\end{equation}
for all $z$ on
$\Omega_3, \Omega_4$, apart from two small arcs
$\Omega_5, \Omega_6$, centred on the ray
$\arg z = L_2(r) + \pi/2 \rho$,
on which  $|A(z)|$ is  bounded above via
(\ref{est3}). By (\ref{x10}), (\ref{x14}) and  (\ref{sigmaest1}),
it may be assumed that the angular measure of each of
$\Omega_5, \Omega_6$ tends to $0$ as $r \to + \infty, r \in F_1$.

Apply the two-constants theorem \cite{Nev} to the subharmonic functions $\log |A + 1/4|$
and $\log |A'|$ on $\Omega$,
and let $z \in \Omega_2$. By rotation,  the harmonic measure
$\omega (z, \Omega_0, \Omega)$ is at least $c(K)> 0$, where $c(K)$ depends only on $K$. On the other hand, since $\Omega$ lies inside the disc
$D(0, K^2r)$ and outside
the disc $D(0, K^{-2} r)$, the monotonicity property of harmonic measure with respect to subdomains yields
$$
\lim_{r \to + \infty, r \in F_1} \sup \{ z \in \Omega_2: \,
\omega(z, \Omega_5 \cup \Omega_6, \Omega) \} = 0.
$$
It follows therefore from (\ref{x4}), (\ref{est1}),  (\ref{est2}), (\ref{est3}) and  (\ref{est4}) that, for
$z \in \Omega_2$, as
$r \to + \infty$ in $F_1$,
\begin{eqnarray*}
\log |A(z)+ 1/4 | &\leq&
-  \, \frac{\pi c(K)}2 \,  L_1(r) (K^{-2} r)^{ \rho } +
3M \log K^2r
+ o( L_1(r) r^\rho )
\\
&\leq&
-  \, \frac{\pi c(K)}4 \, r^{o(1)}  (K^{-2} r)^{ \rho } +
3M \log K^2r
< 0.
\end{eqnarray*}
A similar argument may be applied to $\log |A'|$, using (\ref{est1}) and
 the same upper bounds for $\log |A|$ on $\partial \Omega$, coupled with (\ref{logderest1}).
 This yields,
 provided $r \in F_1$ is large enough,
\begin{equation}
 \label{est5}
 |A(z) | < \frac54 \quad  \hbox{and} \quad |A'(z) | < 1
 \quad
 \hbox{ for $z \in \Omega_2$.}
\end{equation}


The proof of the lemma will now
be completed by estimating
$E(z)$ on $\Omega_2 $. First, (\ref{x4}),
(\ref{x14}) and (\ref{sigmaest1})
yield
\begin{equation}
 \label{est6}
 \log |E(z)|
 \sim \frac{\pi}{\sqrt{2}}
 L_1(r) |z|^\rho
\end{equation}
on $\Omega_2$. Differentiation of the non-linear equation (\ref{bleq}) yields a linear equation
\begin{equation}
E''' + 4AE' + 2A' E = 0.
\label{thirdorder}
\end{equation}
For $j = 0, 1, 2$ write
$$
\lambda = e^{ i (L_2(r) +
3 \pi/4 \rho + \sigma (r) )},
\quad
V_j (t) =
E^{(j)} ( \lambda t ),
\quad V(t) =
|V_0(t)| + |V_1(t)| +
| V_2(t)| ,
$$
so that, by
(\ref{thirdorder}), with the primes denoting $d/dt$,
$$
V_0'(t) = \lambda  V_1 (t),
\quad
V_1'(t) = \lambda  V_2(t),
\quad
V_2' (t)= -
\lambda
(4A(\lambda t) V_1(t) + 2A'(\lambda t) V_0(t)) .
$$
For
$K^{-1}r \leq s \leq K r$,
it follows,
in view of  (\ref{est5}), that
$$
|V_0'(s)| \leq V(s), \quad
|V_1'(s)| \leq V(s), \quad
|V_2'(s)| \leq 5 V(s),
$$
and hence
$$
V(s) \leq v(s) =
V( K^{-1} r ) + \int_{K^{-1} r}^s 10 V(t) \, dt.
$$
Using the standard method from Gronwall's lemma
\cite{Hil1,Hil2} yields
$$
v'(s) = 10 V(s) \leq 10 v(s)
$$
and so, for
$K^{-1}r \leq s \leq K r$,
$$
\log V(s) \leq \log v(s) \leq
\log v(K^{-1} r ) +
10 (s - K^{-1} r) =
\log V(K^{-1} r ) +
10 (s - K^{-1} r).
$$
But this implies, in view of
(\ref{logderest1}) and (\ref{est6}), that
\begin{eqnarray*}
(1-o(1)) \frac\pi{\sqrt{2}} L_1(r) (K r)^\rho &\leq&
\log |E(Kre^{ i (L_2(r) +
3 \pi/4 \rho + \sigma(r) )})|
\\
&\leq& \log V(Kr)
\leq
\log V(K^{-1} r ) +
10 (K - K^{-1}) r \\
&\leq&
(1+o(1)) \frac\pi{\sqrt{2}} L_1(r) (K^{-1} r)^\rho
 + O( \log r) +
10 (K - K^{-1}) r .
\end{eqnarray*}
By (\ref{x4}) and the fact that $r$ may be chosen arbitrarily  large, this forces $\rho =1$.
\hfill$\Box$
\vspace{.1in}

Thus $1$ and $\infty$ are the only deficient values of $E$.
Moreover, the method of Lemma \ref{lemrho1} shows that  if $r > 0$  is large and  the circle $S(0, r)$ does not meet $C_1$, then $S(0, r)$ consists of the following: an arc on which $E(z)$ is close to $1$ and one on which $|E(z)|$ is large,
with
$\log |A(z)| \leq O( \log r )$ on both, by (\ref{bleq});
  two arcs whose  angular measures each tend to $0$ as $r \to + \infty$.
Hence there exists  $N \in \N$ such that if $|z| \geq 1$ lies outside a set of zero logarithmic density, then
$\log |A(z)| \leq N \log |z|$,
apart from on
a set of angular measure $o(1)$. Further, there exists a polynomial $P_1$, of degree at most $N$, such that $A_1 (z) =
(A(z) - P_1(z))/z^{N+1}$ is entire. Applying Lemma \ref{lemtsujiest} to $\log^+ |A_1|$,
which  has order at most $\rho$,  gives a contradiction.
\hfill$\Box$
\vspace{.1in}


\section{Proof of Theorem \ref{thmAsing}}
\label{pfthmAsing}

\begin{prop}
 \label{propAsing}
 Let $A$ be a transcendental entire function, and assume that there exists a simple path $\gamma$ tending to infinity
 such that
 (\ref{de1}) has solutions $u, v$ satisfying
 \begin{equation}
  \label{uvprop}
  u(z) \sim 1, \quad v(z) \sim z, \quad
  W(u, v) \sim 1 ,
 \end{equation}
as $z \to \infty$ on
$\gamma$.
 Let $U = f_2/f_1$ and $E = f_1 f_2$, where $f_1, f_2$ are linearly independent solutions of (\ref{de1}) with $W(f_1, f_2) = 1$.
 Then $U$ does not belong to the Speiser class $\mathcal{S}$.
 In particular, if $A$ has finite order then
 $\lambda(E) = + \infty$.
\end{prop}
Theorem \ref{thmAsing} follows from Proposition \ref{propAsing}, since
\cite[Proposition 7.1]{Lasing2016} guarantees the existence of the required path $\gamma$ and solutions $u, v$.
\\
\textit{Proof of Proposition \ref{propAsing}}.
Let $A$, $E = f_1 f_2$ and $\gamma$ be as in the hypotheses, and assume that $U \in S$,
which will certainly be the case if
$\rho(A) + \lambda (E) < + \infty$,
 by Theorem \ref{thmspeiser}.
Then
$V = v/u$ also belongs to $\mathcal{S}$.
Now (\ref{uvprop}) yields, as $z \to \infty$ on $\gamma$,
$$
V(z) \sim z \to \infty, \quad
V'(z) =
\frac{W(u,v)}{u^2}
\sim 1
, \quad \frac{zV'(z)}{V(z)} \to 1.
$$
But a standard estimate \cite{Ber4,ripstall,sixsmithEL} for functions in the
Eremenko-Lyubich class $\mathcal{B} $, which contains $\mathcal{S}$,  gives
a positive constant $c$ such that, if
$|z|$ and $|V(z)| $ are large enough,
$$
\left| \frac{zV'(z)}{V(z)}
\right| \geq c
\log |V(z)| .
$$
This contradiction proves Proposition \ref{propAsing} and hence Theorem
\ref{thmAsing}.
\hfill$\Box$
\vspace{.1in}

\begin{cor}
 \label{cor9}
 Let $A$ be a transcendental entire function of finite order, and assume that there exists $\theta \in [0, 2 \pi ]$ such that
 \begin{equation}
  \label{thetachoice}
\int_0^{+ \infty}
r |A(r e^{i \theta} )| \, dr < + \infty.
 \end{equation}
Let $E = f_1 f_2$, where $f_1, f_2$ are linearly independent solutions of (\ref{de1}) with $W(f_1, f_2) = 1$.
 Then
 $\lambda(E) = + \infty$.
\end{cor}
In \cite[Theorem 1]{BLL1} it is shown that if $A$ is a transcendental entire function of finite order $\rho$, and almost all $\theta \in [0, 2 \pi]$ are such that either   (\ref{thetachoice}) holds, or $r^{-N} |A(r e^{i \theta} )| \to + \infty $ for every $N \in \N$, or $|A(re^{i\theta})| = O(r^n) $ for some
$n = n(\theta) \geq 0 $ with
$(n+2)/2 < \rho$, then $\lambda(E) = + \infty$. Corollary~\ref{cor9} shows that even a single
occurrence of (\ref{thetachoice}) suffices to force
$\lambda(E) = + \infty$.\\\\
\textit{Proof of Corollary \ref{cor9}.} It suffices to show that (\ref{thetachoice})
implies the existence of  $u, v$  as in
(\ref{uvprop}),
with $\gamma$ the ray $\arg z = \theta$,
and this requires only a standard iterative construction. Assume without loss of generality that $\theta = 0$, and take $X \geq 1$, so large that
$$
\int_X^{+\infty}
r |A(r)| \, dr < \frac12.
$$
Let $u_0 = 0$ and $u_1 = 1$, and assume that $n \in \N$ and
$u_0, \ldots , u_n$ have been constructed with
$$
u_{j}(x) =
1 + \int_{x}^{+\infty}
(x-r) A(r) u_{j-1}(r) \, dr \quad \hbox{and}
\quad
|u_j(x) - u_{j-1}(x)| \leq 2^{1-j}$$
for $x \geq X$, which is evidently true for $n=1$.
Then $|u_n(r)| \leq 2$ for $r \geq X$,
so that
$u_{n+1}$ may be defined and satisfies,
for $x \geq X$,
$$
|u_{n+1}(x) - u_n(x)| \leq
\int_x^{+\infty}
(r-x) |A(r)| |
u_{n}(r) - u_{n-1}(r) | \, dr
\leq 2^{1-n}
\int_x^{+\infty}
r |A(r)| \, dr \leq
2^{-n} .
$$
Hence the series $u(x) =
\sum_{j=1}^\infty (u_j(x) - u_{j-1}(x)) = \lim_{n \to \infty} u_n(x)
$ converges uniformly for $x \geq X$, with
$|u(x)| \leq 2$ there, and satisfies
$$
u(x) =
1 + \int_{x}^{+\infty}
(x-r) A(r) u(r) \, dr ,\quad
u'(x) =
\int_{x}^{+\infty}
 A(r) u(r) \, dr ,\quad u''(x) = -A(x) u(x).
$$
Finally, $u(x) \to 1$ as $x \to + \infty$ and
$v$ may be defined by $(v/u)' = u^{-2}$.
\hfill$\Box$
\vspace{.1in}


The following example shows that
 Corollary \ref{cor9} fails for infinite order. Define an entire function $h$,
a zero-free Bank-Laine function $E$ and a meromorphic function $G$  by
$$
h(z) = \int_1^z \frac{1- e^{-t}}t \, dt
+ \int_1^{+ \infty } \,
\frac{dt}{t e^t} ,
\quad E = e^h , \quad
\frac1{E(z)} =
e^{-h(z)} =
 \frac1z + G (z)
.
$$
Let $D$ be the sectorial domain given by $1/2 < |z| < + \infty$, $| \arg z | < \pi /4$, and denote by $\varepsilon(z)$ any term with the property that
$z^n \varepsilon(z) \to 0$ uniformly as
$z \to \infty$ in $D$, for every $n \in \N$. Then
$$
h(z) = \log z +
\int_z^{+ \infty}
\,
\frac{dt}{t e^t} =
\log z + \varepsilon (z), \quad
h'(z) = \frac{1- e^{-z}}z = \frac1z + \varepsilon(z),
\quad h''(z) = - \,
\frac{1}{z^2}
+ \varepsilon(z),
$$
in which the integral is eventually along the positive real axis.
Now define an entire function $A$ by (\ref{bleq}). Then $E$ and $G$ satisfy
$$
E(z) = z + \varepsilon (z), \quad
G(z) = \varepsilon (z) ,
$$
and (\ref{bleq}) delivers
(\ref{thetachoice}) with $\theta = 0$, via
$$
- 4A(z) = h'(z)^2
+ 2 h''(z) + e^{-2h(z)} = - \,
\frac{1}{z^2}
+ \varepsilon(z)
+
\frac{1}{z^2} (1 +
 \varepsilon(z)) =
\varepsilon (z).
$$
Finally, define
$ g_1, g_2$ on $D$ by
\begin{eqnarray*}
g_j(z) &=& E(z)^{1/2}
\exp \left( \frac{(-1)^j}2
\left( \log z - \int_z^{+ \infty}
\, G(t) \, dt
\right) \right)
\\
&=&
\exp \left( \frac{h(z)}2 +  \frac{(-1)^j}2
\left( \log z - \int_z^{+ \infty}
\, G(t) \, dt
\right) \right)
,
\\
\frac{g_j'(z)}{g_j(z)} &=&
\frac12 \left(
\frac{E'(z)}{E(z)}
+ (-1)^j \left( \frac1z + G(z) \right) \right)
=
\frac12 \left(
\frac{E'(z)}{E(z)}
+ \frac{(-1)^j }{E(z)}  \right).
\end{eqnarray*}
Straightforward computation then shows that $g_1, g_2$ solve (\ref{de1}) and
$W(g_1, g_2) = 1$,
while
$$
g_1(z) = \exp( \varepsilon (z) ) =
1 + \varepsilon(z), \quad g_2(z) =
\exp( \log z + \varepsilon (z) )
= z + \varepsilon (z).
$$
\hfill$\Box$
\vspace{.1in}




\section{Proof of Theorem \ref{thmdefA}}
\label{pfthmdefA}

The first step is to establish the following proposition.

 \begin{prop}
 \label{propBsing}
  Let $U$ be a
  locally univalent transcendental meromorphic function in the plane, and assume that $U \in \mathcal{S}$ and
  that the Schwarzian derivative $S(U)$ is transcendental of finite order.
  Then $\delta (0, S(U)) = 0$.
 \end{prop}
 \textit{Proof.}
 Let $n_0 < + \infty$ be the number of  singular values of
 $U^{-1}$ in
 $\C \cup \{ \infty \}$, and define $A$ by $S(U) = 2A$,
 so that $U$ is the quotient of linearly independent solutions of (\ref{de1}).
 It is clear that $\delta (0, A)) = \delta (0, S(U))$, so assume  that $\delta(0, A) =
 \delta > 0$.
  The first step is to determine  positive real numbers $\delta_1, K_1$, as well as  sequences $r_k \to + \infty$, $\theta_k \in \R$,
  such that
  \begin{equation}
   \label{T(r)est}
   T(8 r_k, A) < K_1 T(r_k, A)
  \end{equation}
and
 \begin{equation}
 \log |A(r_k e^{i \theta} )| < - \, \frac{\delta}4 \,  T(r_k, A)
 \quad \hbox{for} \quad
 \theta_k \leq \theta \leq \theta_k + \delta_1.
 \label{Aest0}
 \end{equation}
To prove (\ref{T(r)est}) and (\ref{Aest0}), first
write $A(z) = \alpha_1 z^{n_1} A_1(z)$, with $\alpha_1 \in \C $ and $n_1 \in \Z$, so that $A_1(0) = 1$, and observe that
$\delta(0, A_1) =
\delta > 0$. Then choose $q_k \to + \infty$ with
$T(16 q_k, A) < K_1 T(q_k, A)$, which is possible for large enough $K_1$ since $A$ has finite order. Next, apply Lemma~\ref{lemfuchssmallarcs} to $g = A_1$, with $R = 16q_k$, $\eta_1 = 1/2$,  and with $ \eta_2$  sufficiently small that there exists $r = r_k \in [q_k, 2 q_k] \cap F_R$, which implies
(\ref{T(r)est}).
Finally, take $G_r =
[\theta_k, \theta_k + \delta_1] \subseteq [0, 2 \pi]$,  with $\theta_k$ chosen so that
$\log |A_1(z)| <
- \delta T(r_k, A_1)/2$ at one end-point of the arc
$\{ r e^{i \theta} :
\, \theta \in G_r \}$.
If
$m = \delta_1 $ is chosen small enough, (\ref{smallarcs})  delivers (\ref{Aest0}). It may be assumed that $(\theta_k)$ converges and, after a rotation of the independent variable, that
 \begin{equation}
  \label{Aest1}
  \log |A(r_k e^{i \theta} )| < - \, \frac{\delta}4 \, T(r_k, A)
 \quad \hbox{for} \quad
 - \delta_2 \leq  \theta   \leq
 \delta_2 , \quad
 \delta_2 = \frac{\delta_1}4.
 \end{equation}
The idea of the proof is to show that $2A$ can be written as the Schwarzian derivative of a function $v_k/u_k$ in the Speiser class which is close to the identity on the arc in (\ref{Aest1}).
After normalising these quotients $v_k/u_k$,
continuation of the  inverse functions
into  suitable regions will deliver  a sequence of circles $|z| = s_k r_k$, on nearly all of which $A(z)$ is small, which will then lead to a contradiction via
(\ref{T(r)est}) and Lemma \ref{EFlemma}.
The  proof will use $c, K$  to denote positive constants, not necessarily the same at each occurrence, but always
  independent of $r_k$, with $c$ small and $K$ large.


\begin{lem}
 For large $k$ there exist solutions
 $u_k, v_k$ of (\ref{de1}) such that $W(u_k, v_k) = 1$ and
 \begin{equation}
  \label{Aest2a}
  v_k(z) =
 z \left( 1 + \varepsilon_{1,k}(z) \right),
 \quad
 u_k(z) =
  1 + \varepsilon_{2,k}(z),
 \end{equation}
 in which
 \begin{equation}
  \label{Aest2}
 | \varepsilon_{j,k}(r_ke^{i\theta}) |
 < \exp( - c\delta
 T(r_k, A) )
 \quad \hbox{for} \quad
  - \delta_2 \leq  \theta   \leq
 \delta_2 .
 \end{equation}
\end{lem}
\textit{Proof.}
This is standard.
Let $k$ be large, denote the arc
$\{ r_k e^{i \theta } :
\,  - \delta_2 \leq  \theta   \leq
 \delta_2 \}$ by
$B_k$ and let $z_k$ be one of its end-points. If $u$ is any solution of (\ref{de1}) write
$$
v(z) = \int_{z_k}^z
(t-z) A(t) u(t) \, dt
,\quad
v'(z) = - \int_{z_k}^z
 A(t) u(t) \, dt,
 \quad v''(z) = -A(z) u(z) = u''(z) .
$$
Since $v(z_k) = v'(z_k) = 0$ this yields
\begin{equation}
 \label{uident}
u(z) = u(z_k) +
u'(z_k) (z-z_k) + v(z) = u(z_k) +
u'(z_k) (z-z_k) +
\int_{z_k}^z
(t-z) A(t) u(t) \, dt .
\end{equation}
If $\max \{ |u(z_k)| , |u'(z_k)| \} \leq  r_k$ then
 $|u(z)| < r_k^3$ on $B_k$.
To see this,
suppose that $z_k'$ is the nearest point to $z_k$ on $B_k$ with
$|u(z_k')| \geq r_k^3$; then the integral in (\ref{uident}) is small at $z = z_k'$,
by (\ref{Aest1}), giving $|u(z_k')| \leq 4r_k^2$, a contradiction.
Now choose solutions $v_k, u_k$ of (\ref{de1})
such that
$$
u_k(z_k) = 1, \quad u_k'(z_k) = 0, \quad
v_k(z_k) = z_k, \quad
v_k'(z_k) = 1,
$$
and estimate the integral in (\ref{uident}) using (\ref{Aest1}) again.
\hfill$\Box$
\vspace{.1in}


There exists a M\"obius transformation $S_k$ such that
$v_k/u_k = S_k \circ U$. Write
\begin{equation}
 \label{Aest3}
V_k( \zeta ) =
\frac{S_k (U( r_k \zeta))}{r_k}
= \frac{v_k(r_k \zeta)}{r_k u_k(r_k \zeta)} =
\zeta (1 + \rho_k(\zeta)) ,
\quad V_k'(\zeta) =
\frac1{u_k(r_k \zeta)^2} =
1 + \sigma_k(\zeta),
\end{equation}
and denote by $P_k$ the Schwarzian
$P_k = S(V_k)$ of $V_k$.
Then (\ref{Aest1}),
(\ref{Aest2a}) and (\ref{Aest2})
imply that
\begin{equation}
 \label{Aest4}
 | \rho_k(\zeta) |
 + | \sigma_k(\zeta)|
 + | P_k(\zeta) |
 < \exp( - c \delta
 T(r_k, A) )
 \quad \hbox{for $\zeta \in \Omega_k =
 \{ e^{i \theta} : \,
  - \delta_2 \leq  \theta   \leq
 \delta_2 \}$.}
\end{equation}


The function $V_k$ is transcendental and meromorphic in the plane and, since $S_k$ is M\"obius, the inverse of $V_k$ has  $n_0$ singular values. Hence there exist
$\phi_k \in [ - \delta_2/2, \delta_2/2]$  and $s_k \in [2, 3]$ such that
$V_k$ has no asymptotic or critical values $w$ in the simply connected
domain $D_k$ given by
$$
D_k = \{ w \in \C \setminus \{ 0 \} : \, | \arg w - \phi_k |
<  \delta_3 \} \cup
\{ w \in \C : \, ||w| - s_k| <  \delta_3,
 \, | \arg w - \phi_k |
< \pi \},
$$
in which $\delta_3 > 0$ is small compared to $\delta_2$, but independent of $k$.
In particular, $V_k$ is univalent on any component
of the pre-image of  $D_k$.

Let $\delta_4 > 0$ be small compared to $\delta_3$.
For large $k$, by (\ref{Aest3}) and (\ref{Aest4}), there exists $\zeta_k$ on the arc $\Omega_k$
 such that
\begin{equation}
 \label{zetakdefA}
V_k(\zeta_k) = w_k,
\quad
\frac45 < |w_k| < \frac54, \quad
\arg w_k = \phi_k .
\end{equation}
Further, the disc
$D(w_k, 4 \delta_4)$
lies in $D_k$ and
there exists
a sub-arc
$\Omega_k'$ of $\Omega_k$  mapped univalently by $V_k$ onto a simple path $\Omega_k''$
lying in $|w-w_k| \leq 2 \delta_4$ and joining  $w_k$ to the circle $S(w_k,  2 \delta_4)$.
Let $W_k$ be that branch of $V_k^{-1}$ which maps $w_k$ to $\zeta_k$: then $W_k$ extends univalently to $D_k$ and $W_k(V_k(\zeta)) = \zeta$ on $\Omega_k'$.
Moreover, $W_k$ and its Schwarzian $Q_k = S(W_k)$ satisfy, by (\ref{Aest3}), (\ref{Aest4}) and the composition formula
\begin{equation*}
 0 = S(W_k \circ V_k) =
Q_k(w) V_k'(\zeta)^2
+ P_k(\zeta) ,
\end{equation*}
the estimate
\begin{equation}
 \label{Aest5}
 |W_k'(w) - 1 | +
|Q_k(w)| <
\exp( - c \delta
 T(r_k, A))
 \quad \hbox{on
 $\Omega_k''$.}
\end{equation}

\begin{lem}
 If $k$ is large enough then
 $W_k$ and $Q_k$
 satisfy
 \begin{equation}
 \label{Aest5a}
 |W_k'(w) - 1 | +
|Q_k(w)| <
\exp( - c \delta
 T(r_k, A))
 \quad
 \hbox{on the disc $D(w_k, \delta_4)$.}
\end{equation}
\end{lem}
\textit{Proof.} First,
$W_k$ is univalent on
 $D(w_k, 4 \delta_4)  \subseteq D_k$, which delivers an estimate
\begin{equation}
 | W_k'(w)| +
 |Q_k(w)| < K
\label{koebenehari}
\end{equation}
on $D(w_k, 2 \delta_4)$:
for $W_k'$ this follows from
Koebe's distortion theorem \cite{Hay9} and
(\ref{Aest5}), while for $Q_k$ it suffices to use Nehari's univalence criterion \cite{nehari}.
The estimate (\ref{Aest5}) holds on the simple path $ \Omega_k''$ from $w_k$ to
$S(w_k, 2 \delta_4)$, and for
$w \in D(w_k, \delta_4) \setminus \Omega_k''$ the harmonic measure of $\Omega_k''$ is at least~$c$:
thus the two constants theorem
\cite{Nev} applied to the subharmonic function
$\log |W_k'-1|$
gives
$$
\log |W_k'(w)-1| < K - c \delta T(r_k, A)
$$
on $D(w_k, \delta_4)$. Applying the same argument to
$\log |Q_k|$
yields (\ref{Aest5a}).
\hfill$\Box$
\vspace{.1in}

\begin{lem}
\label{lemkoebe1}
Let $\delta_5 > 0$ be small compared to $\delta_4$.
 If $k$ is large enough then
 $W_k$ and $Q_k$
 satisfy an estimate
 (\ref{koebenehari})
 on the sub-domain $D_k' = D_{1,k}' \cup D_{2,k}'$ of
 $D_k$, where
 \begin{eqnarray*}
D_{1,k}' &=&
\{ w \in \C \setminus \{ 0 \} : \, |w_k|/32 < |w| < 32|w_k|, \, | \arg w - \phi_k |
<  2 \delta_5 \} ,\\
D_{2,k}' &=&
\{ w \in \C : \, ||w| - s_k| <  2 \delta_5,
\,
| \arg w - \phi_k | < \pi - \delta_5 \}.
\end{eqnarray*}
\end{lem}
\textit{Proof.} The estimate for $Q_k$ again follows from Nehari's criterion
\cite{nehari}, since  for each
$w \in D_k'$ the function
$W_k$ is univalent on a disc $D(w, \delta_6) \subseteq D_k$
of  fixed radius $\delta_6$.
To estimate $W_k'(w)$, take a small $\delta_7 > 0$,
start from
(\ref{Aest5}) at the point $w_k \in \Omega_k''$ and apply Koebe's distortion theorem \cite{Hay9} repeatedly on a chain of discs
$D(w_{k,j}, \delta_{7}) $, $0 \leq j \leq K$, the centres being chosen so that
$w_{k,0} = w_k$ and the last $w_{k,j} $ is $w$, while
$w_{k,j+1} \in D(w_{k,j}, \delta_{7}/2)$.
\hfill$\Box$
\vspace{.1in}

\begin{lem}
\label{lemkoebe}
 If $k$ is large enough then
 $W_k$ and $Q_k$
 satisfy
 \begin{equation}
 \label{Aest6}
 |W_k'(w) - 1 | +
|Q_k(w)| <
\exp( - c \delta
 T(r_k, A))
\end{equation}
 for all $w$ in the sub-domain $D_k'' =
 D_{1,k}'' \cup D_{2,k}''$ of
 $D_k'$,
 where
 \begin{eqnarray*}
D_{1,k}'' &=&
\{ w \in \C \setminus \{ 0 \} : \, |w_k|/16 < |w| < 16 |w_k|, \, | \arg w - \phi_k |
<  \delta_5 \} ,\\
D_{2,k}'' &=&
\{ w \in \C : \, ||w| - s_k| <   \delta_5,
\,
| \arg w - \phi_k | < \pi - 2 \delta_5 \}.
\end{eqnarray*}
\end{lem}
\textit{Proof.}
The subharmonic functions
$\log |W_k'-1|$ and
$\log|Q_k|$ are uniformly bounded above on $D_k'$, by Lemma \ref{lemkoebe1}, and the stronger estimate (\ref{Aest5a}) holds on $D(w_k, \delta_4)$, and so
on a subset of $\partial D_{1,k}'$ whose harmonic measure with respect to
$D_{1,k}'$ is at least $c$ when evaluated at any
$w \in D_{1,k}''$.
Thus the two constants theorem
delivers (\ref{Aest6})
on
$D_{1,k}''$.
This estimate may then be extended
in the same way
to
$D_{2,k}''$,
possibly with a smaller constant $c$.
\hfill$\Box$
\vspace{.1in}

\begin{lem}
\label{lemWkimage}
 If $k$ is large enough then
$$
W_k(w) =
w ( 1 + \tau_k(w)) ,
\quad
| \tau_k(w) | <
\exp( - c \delta
 T(r_k, A))
 \quad \hbox{on $D_k''$}.
 $$
\end{lem}
\textit{Proof.} Write
$\tau(w)$ for any term of modulus less than
$\exp( - c \delta
 T(r_k, A))$. Then (\ref{Aest3}),
 (\ref{Aest4}),
 (\ref{zetakdefA})
 and integration of
(\ref{Aest6}) yield
\begin{eqnarray*}
 W_k(w) &=& W_k(w_k) +
 (w-w_k)(1 + \tau(w))
 = \zeta_k +
 (w-w_k)(1 + \tau(w))\\
  &=& w_k(1 + \tau(w)) +
 (w-w_k)(1 + \tau(w))
 = w(1 + \tau(w)).
\end{eqnarray*}
\hfill$\Box$
\vspace{.1in}

\begin{lem}
\label{lemschwarzianest}
 If $k$ is large enough then the image of $D_k''$ under $W_k$ contains the arc
 $\Lambda_k$ given by
 $| \zeta| = s_k$,
$| \arg z - \phi_k |
\leq \pi - 3 \delta_5$,
and $P_k = S(V_k)$ satisfies
$$
| P_k(\zeta) |
< \exp( - c \delta
 T(r_k, A)) \quad\hbox{on $\Lambda_k$.}
$$
\end{lem}
\textit{Proof.} This holds by
(\ref{Aest6}),
Lemma \ref{lemWkimage} and the fact that
$V_k \circ W_k$ is the identity on $D_k$.
\hfill$\Box$
\vspace{.1in}

It now follows from Lemma \ref{lemschwarzianest} that
$$
| A(z) |
< \exp( - c \delta
 T(r_k, A)) \quad\hbox{for $|z| = s_k r_k$, $| \arg z - \phi_k| \leq
 \pi - 3 \delta_5 $.}
$$
Since $2 \leq s_k \leq 3$ this implies that, by (\ref{T(r)est})
and Lemma \ref{EFlemma},
\begin{eqnarray*}
 T(r_k s_k, A) =
 m(r_k s_k, A)
 &\leq&
 \frac{528 r_k \delta_5}{8r_k -s_k r_k} \left( 1 + \log
 \frac1{6 \delta_5} \right) T(8r_k, A)\\
 &\leq& 132 K_1 \delta_5 \left( 1 + \log
 \frac1{6 \delta_5} \right)T(r_k, A)
\\
 &\leq& 132 K_1 \delta_5 \left( 1 + \log
 \frac1{6 \delta_5} \right)T(r_k s_k, A).
\end{eqnarray*}
But $K_1$ is independent of $\delta_5$, which may in turn be chosen arbitrarily small, and this contradiction proves Proposition \ref{propBsing}.
\hfill$\Box$
\vspace{.1in}

To prove
Theorem \ref{thmdefA},
suppose that $A$ is a transcendental entire function of finite order and that
 $\lambda (f_1f_2 )
 < + \infty $, where $f_1, f_2$ are linearly independent solutions of (\ref{de1}). Then $E = f_1 f_2$ has finite order
 and $U = f_2 /f_1$ belongs to the Speiser class
 $\mathcal{S}$,  by Theorems \ref{thmA} and \ref{thmspeiser}.
 As already observed,
 $\delta (0, A) =
 \delta (0, S(U))$,
 where $S(U) = 2A$ is the Schwarzian of $U$ as in (\ref{schwarzian}),
 and Proposition
 \ref{propBsing}
 forces
 $\delta (0, A)  = 0$.
\hfill$\Box$
\vspace{.1in}


\subsection{Proof of Corollary \ref{corA}} By Theorem
\ref{thmdefA} it suffices to prove that $\delta(0, A) > 0$. To this end, write
$C = f'/f$ and
$$
- A = C' + C^2, \quad
- \, \frac1C =
\frac1{A} \left( \frac{C'}C + C \right),
$$
so that standard properties of $\log^+$ deliver
\begin{equation}
\log^+ \frac1{|C(z)|}
\leq
\log^+ \frac1{|A(z)|} +
\log^+ \frac{|C'(z)|}{|C(z)|} + \log 2.
\label{1/Cestimate}
\end{equation}
It follows that, outside a set of finite measure,
$$
\frac{ \delta (0, C)}
2 \, T(r, C) \leq
m(r, 1/C) \leq
m(r, 1/A) + m(r, C'/C) + \log 2 \leq
T(r, A) + o( T(r, C)).
$$
Thus $C$ has finite order, by the  monotonicity of $T(r, C)$. Using (\ref{1/Cestimate}) again  yields, for all large $r$,
\begin{eqnarray*}
T(r, A) &\leq&
O( T(r, C)) + O( \log r)  \leq
O( m(r, 1/C)) \\
&\leq&
O( m(r, 1/A) ) + O( \log r )
\leq
O( m(r, 1/A) ) + o( T(r, A) ).
\end{eqnarray*}
\hfill$\Box$
\vspace{.1in}

\subsection{Proof of Corollary \ref{corB}} Suppose that
$\lambda(E) < + \infty$. Then dividing (\ref{specialde}) by $E$ and applying the lemma of the logarithmic derivative delivers
$$
m(r,1/E) \leq m(r, B) + o( T(r, E))
$$
outside a set of finite measure, and hence
$\rho(E)  < + \infty$.
Now $E$ is a Bank-Laine function, and so there exists an entire function $A$ such that $E = f_1 f_2$, with the $f_j$ normalised linearly independent solutions of (\ref{de1}), and $A$ has finite order by (\ref{bleq}).
With $U = f_2/f_1$ this gives
$$
\frac{E'}{E} =
\frac{f_1'}{f_1} +
\frac{f_2'}{f_2}, \quad
\frac1E = \frac{U'}{U} =
\frac{f_2'}{f_2} -
\frac{f_1'}{f_1} ,
\quad B = \frac{E'-1}E =
2 \, \frac{f_1'}{f_1},
$$
and so
$\delta (0, f_1'/f_1) > 0$. Since $\lambda (f_1 f_2 ) =
\lambda (E) < + \infty$, this contradicts Corollary \ref{corA}.
\hfill$\Box$

\noindent
J.K. Langley, Emeritus Professor,\\
School of Mathematical Sciences, University of Nottingham, NG7 2RD, UK.\\
james.langley@nottingham.ac.uk

\end{document}